\documentclass[12pt]{article}
\usepackage{latexsym}
\usepackage{amssymb}
\begin{document}
\newtheorem{theorem}{Theorem}
\newtheorem{lemma}{Lemma}
\newtheorem{proposition}{Proposition}
\newtheorem{corollary}{Corollary}
\newtheorem{definition}{Definition}
\newtheorem{question}{Question}
\newtheorem{conjecture}{Conjecture}
\newcommand{\F}{\ensuremath{\mathbb F}}
\newcommand{\N}{\mathcal N}
\newcommand{\R}{\mathcal R}

\title{Restrictions on sets of conjugacy class sizes in arithmetic progressions}
\author{Alan R. Camina \& Rachel D. Camina}
\date {}
\maketitle

\begin{abstract}
We continue the investigation, that began in \cite{bianchi} and \cite{glasby}, into finite groups whose set of nontrivial conjugacy class sizes form an arithmetic progression.
Let $G$ be a finite group and denote the set of conjugacy class sizes of $G$ by ${\rm cs}(G)$. Finite groups satisfying 
${\rm cs}(G) = \{1,2,4,6\}$ and $\{1,2,4,6,8\}$ are classified in \cite{glasby} and \cite{bianchi}, respectively, we demonstrate these examples are rather special by proving the following. There exists a finite group $G$ such that 
${\rm cs}(G) = \{1, 2^{\alpha}, 2^{\alpha+1}, 2^{\alpha}3 \}$ if and only if $\alpha =1$. Furthermore, there exists a finite
group $G$ such that
 ${\rm cs}(G) = \{1, 2^{\alpha}, 2^{\alpha +1}, 2^{\alpha}3, 2^{\alpha +2}\}$ and $\alpha$ is odd if and only if $\alpha=1$.

\end{abstract}

\section{Introduction}

Let $G$ be a finite group and ${\rm cs}(G)$ the set of conjugacy class sizes of elements of $G$. That is, 
${\rm cs}(G) = \{ |x^G| : x \in G \}$ where $x^G$ denotes the conjugacy class of $x$ in $G$. 
The question as to which sets of natural numbers can occur as ${\rm cs}(G)$ for some finite group $G$, and also
the relation between ${\rm cs}(G)$ and the algebraic structure of $G$, are  questions 
that have long interested mathematicians, see \cite{survey} for an overview of results.
 For example in 1953  It{\^o} \cite{ito} proved that if ${\rm cs}(G) = \{1, n\}$ then $n=p^a$ for some prime $p$ and $G \cong P \times A$ where $P$ is a $p$-group and $A$ is an abelian $p'$-group, more recently Ishikawa proved that such a group
has nilpotency class at most 3 \cite{ishi}.
It is worth noting that if $A$ is an abelian group then ${\rm cs}(G \times A) = {\rm cs}(G)$ and thus, when working with 
${\rm cs}(G)$, we are always working modulo abelian direct factors.
In \cite{cossey} the authors prove that for a given prime $p$ if $S$ is any set
of $p$-powers containing 1 then $S = {\rm cs}(P)$ for some finite $p$-group $P$ of nilpotency class 2. However, in general, there are many restrictions
on the sets that can occur as  sets of conjugacy class sizes. In particular,  in \cite{casolo}
 the authors prove that if the largest two conjugacy
class sizes $m$ and $n$ are coprime, then ${\rm cs}(G) = \{1, m, n\}$ and $G$ is quasi-Frobenius (that is
$G/Z(G)$ is a Frobenius group). 

 Recently there has been interest in when ${\rm cs}(G)$  or 
${\rm cs}^*(G) = {\rm cs}(G) \setminus \{1\}$ can be an arithmetic progression and, when this is so, the algebraic implications on the group
have been investigated. 
That is, we ask when can ${\rm cs}(G)$ or ${\rm cs}^*(G)$ be of the form $\{a, a+d, a+2d, \ldots, a+kd \}$ for natural numbers $a, d$
and $k$?
A simple argument shows that if ${\rm cs}(G)$
is an arithmetic progression then we are in the situation described above, that is the largest two conjugacy class sizes are
coprime and thus ${\rm cs}(G) = \{1, m, n\}.$  This situation has been analysed further in \cite{glasby} where the authors prove
that if ${\rm cs}(G)$ is an arithmetic progression (with at least 3 terms) then necessarily ${\rm cs}(G) = \{1,2,3\}$. The authors
classify all such $G$. 

A less restrictive condition is to insist ${\rm cs}^*(G)$ is an arithmetic progression. 
We call $|{\rm cs}^*(G)|$ the conjugate rank of $G$ and note that groups with conjugate rank 2 have been classified 
\cite{jabara}. We suspect it is quite rare for ${\rm cs}^*(G)$ to be an
arithmetic progression of size at least 3, simple arguments (see the next section) show that such a group must have a centre but cannot be nilpotent.
 In \cite{glasby} the authors classify all groups which
satisfy ${\rm cs}^*(G) =\{2,4,6\}$ and
groups which satisfy 
${\rm cs}^*(G) = \{2,4,6,8\}$ are classified in \cite{bianchi}.  It is interesting to note that if $G$ is a group with order divisible by
only two primes, conjugate rank at least 4  and ${\rm cs}^*(G)$ an arithmetic progression then 
${\rm cs}^*(G) = \{m, 2m, 3m, 4m\}$ for some $m= 2^{\alpha}3^{\beta}$ and 
 natural numbers $\alpha$ and $\beta$ (this follows from the analysis of primitive arithmetic progressions in
\cite[Lemma 5]{glasby}). Thus we are led to considering this scenario, 
we restrict to the case when $\beta =0$ and also consider the conjugate rank 3 case. We
 prove the following theorem, indicating that the examples
above are rather special. Our proofs are quite different to those used in \cite{glasby} and \cite{bianchi}.\\[2ex]
{\bf Theorem.} {\it There exists a finite group $G$ with \\
(i) ${\rm cs}(G) =\{1, 2^{\alpha}, 2^{\alpha +1}, 2^{\alpha}3\}$ if and only if
 $\alpha=1$.\\
(ii) ${\rm cs}(G) = \{1, 2^{\alpha}, 2^{\alpha +1}, 2^{\alpha}3, 2^{\alpha +2}\}$  and $\alpha$ odd if and only if $\alpha =1$}.\\[1ex]

\noindent We note that, for $\alpha=1$, in $(i)$ the Sylow 3-subgroup of $G$ is normal, whereas in $(ii)$ the Sylow 2-subgroup of
$G$ is normal. 
The $\alpha$ even, conjugate rank 4 case remains open. If such a $G$ were to exist it would follow, from the analysis in this paper, that neither the Sylow 2 or 3-subgroup would be normal and moreover $G/F(G)$ would be isomorphic to $S_3$, where $F(G)$ denotes the Fitting subgroup of $G$. Thus we are led to the following (increasingly general) questions.

\begin{question}  Suppose $\alpha$ is even. Does there exist a finite group $G$ satisfying ${\rm cs}^*(G) = \{ 2^{\alpha}, 2^{\alpha+1}, 2^{\alpha}3, 2^{\alpha+2} \}$? \end{question}

\begin{question} For which $m = 2^{\alpha}3^{\beta}$ does these exist a finite group $G$ satisfying 
${\rm cs}^*(G) = \{ m, 2m, 3m, 4m\}$? \end{question} 

\begin{question} For which natural numbers $a, k$ and $d$ does there exist a finite group $G$ satisfying 
${\rm cs}^*(G) = \{a, a+d, a+2d, \ldots, a+ kd\}$?
\end{question}

\noindent {\bf Notation and terminology} used throughout the paper are mostly standard. For example for $G$ a group and $x \in G$ we denote
the centraliser of $x$ in $G$ by $C_G(x)$, 
 $Z(G)$ denotes the centre of  $G$, 
$\Phi(G)$ is the Frattini subgroup of $G$ and $O_2(G)$ denotes the $2$-core of $G$, that is, the largest normal 2-subgroup in $G$. However we also use 
the terminology of It{\^o} and call $|x^G|$ the {\it index of $x$  in G} or simply the {\it index of x} when $G$ is clear. Similarly 
${\rm cs}(G)$ is called
the {\it set of indices of $G$}.
 All groups considered are finite.\\[2ex]
\textbf{Acknowledgements} We thank Mariagrazia Bianchi for introducing us to this interesting problem. The second author would like to thank the Isaac Newton Trust for supporting her sabbatical when much of this research took place.

\section{Preliminaries \& Discussion}

Thoughout this paper we will assume well-known properties of indices. Namely that if 
$x$ and $y$ are commuting elements of coprime order then $C_G(xy) = C_G(x) \cap C_G(y)$, that if
$N$ is a normal subgroup of $G$ then both $|x^N|$ and $|(xN)^{G/N}|$ divide $|x^G|$ and finally
if all indices of a
group $G$ are coprime to a prime $p$ then the Sylow $p$-subgroup of $G$ is an abelian direct factor \cite{alan}.   
We start with a definition.
\begin{definition} We say a group $G$  has property ${\mathcal AP}$, or is an ${\mathcal AP}$ group, if $G$ has 
conjugate rank at least 3 and
${\rm cs}^*(G) = \{a, a+d, a+2d, \ldots, a+kd\}$ for some natural numbers $a,d$ and $k$.
\end{definition} 

The following lemma shows that if $G$ is an ${\mathcal AP}$ group then there exists a prime $p$ which divides all
nontrivial indices of elements of $G$.

\begin{lemma}\label{coprime} Suppose ${\rm cs}^*(G) = \{a, a+d, a+2d, \ldots, a+kd\}$ with $a$ and $d$ coprime, then
$|{\rm cs}^*(G)| \leq 2$.\end{lemma} 
{\bf Proof.} Note $\gcd(a+(k-1)d, a+kd) = \gcd (a, d)=1$. The result follows from \cite{casolo}. $\Box$

\begin{lemma} Suppose $p$ is a prime and $p^{\alpha}$ divides the index of all noncentral elements of a finite group $G$.
Let $P$ be 
a Sylow $p$-subgroup of $G$ then $Z(P) = P \cap Z(G)$ and $p^{\alpha}$ divides $|Z(P)|$.
\end{lemma}
{\bf Proof.} By the class equation $p^{\alpha}$ divides $|Z(G)|$. Also, 
any element in the centre of a Sylow $p$-subgroup must be central and thus
$Z(P) = P \cap Z(G)$. The result follows. $\Box$\\

The previous two lemmas show that if $G$ is an ${\mathcal AP}$ group then $G$ has a nontrivial centre.
The next lemma shows that nilpotent groups do not have property ${\mathcal AP}$.

\begin{lemma} Suppose ${\rm cs}^*(G) =\{a, a+d, a+2d, \ldots, a +kd\}$ and $|{\rm cs}^*(G)| \geq 3$. Then $G$ is not nilpotent.
\end{lemma}
{\bf Proof.} Suppose $G$ is nilpotent.  
As $|{\rm cs}^*(G)| \geq 3$ we know by Lemma \ref{coprime} that there exists a prime $p$ such that $p$ divides all
nontrivial indices of $G$, this forces $G$ to be a $p$-group. Thus ${\rm cs}^*(G) = \{p^{a_1}, p^{a_2}, \ldots,
p^{a_k}\}$ for some positive integers $a_1 < a_2 < \cdots < a_k$. But then $p^{a_3} = p^{a_2} + d$ and so $p^{a_2}$ divides $d$. However
also $p^{a_2} = p^{a_1} + d$ which leads to $p^{a_2}$ dividing $p^{a_1}$, a contradiction.$\Box$\\

A group is called an $F$-group if given any noncentral elements $x$ and $y$ then $C_G(x) \not< C_G(y)$. It follows that
a sufficient condition to be an $F$-group is to require that
given any noncentral elements $x$ and $y$ then 
$|y^G|$ does not divide $|x^G|$. The $F$-groups have
been classified \cite{rebmann}, a check of the list shows that no $F$-group has property ${\mathcal AP}$.  Thus any example of an ${\mathcal AP}$  group will have
to include some divisibility of indices. However even with divisibility we run into obstacles. 
Applying the main result in \cite{dolfi} shows that many arithmetic progressions cannot occur as sets
of indices. For example it follows that there does not exist a group $G$ with ${\rm cs}^*(G) = \{2,6,10,14\}$ or
more generally with ${\rm cs}^*(G) = \{a, 3a, 5a, 7a\}$ and $a$ coprime to 105.

We now focus on proving our main result. The following lemmas will be useful.

\begin{lemma}\label{central} If $p$ is the highest power of $p$ which divides any index of $G$ 
then $\Phi(P)$ is central in $G$ for any Sylow $p$-subgroup $P$ of $G$.
\end{lemma}
{\bf Proof.} Let $x \in G$ then there exists a Sylow $p$-subgroup
$\hat{P}$ and a subgroup $P_0$ of index $p$ in $\hat{P}$ such that
$C_G(x) \geq P_0$. There exists $g \in G$ such that $\hat{P}^g =P$
and thus $C_G(x^g) \geq P_0^g \geq \Phi(P)$. So $x^g \in C_G(\Phi(P))$.
So, we have shown that any element of $G$ is conjugate to an element in
$C_G(\Phi(P))$ and thus $C_G(\Phi(P))=G$ by a result of Burnside \cite[\S26]{burnside}.$\Box$\\

The following is adapted from \cite[Proposition 1]{camina}.

\begin{lemma}\label{centre} Suppose $P$ is an abelian Sylow $p$-subgroup and $Z$ is a central $p$-subgroup of a group $G$. Then
 ${\rm cs}(G) = {\rm cs}(G/Z)$.
\end{lemma}
{\bf Proof.} Let $K$ be a central $p$-subgroup of order $p$ and denote $G/K$ by ${\bar G}$. Suppose there exists an element
$x \in G$ such that $|x^G| \neq |\bar{x}^{\bar G}|$. Then $x^g = xz$ for some $g \in G$ where $K = \langle z \rangle$.  As
the order of $x$ equals the order of $x^g=xz$  it follows that $p$, the order of $z$, divides
the order of $x$. Write $x=x_{p} x_{p'}$ as a product of its $p$- and $p'$-parts.  By considering $x^m$, where $m$ is the
order of $x_{p'}$, we can assume $x$ is a $p$-element and thus has $p'$-index. Now $x^{g^p} = xz^p = x$ and
so $g^p \in C_G(x)$. Writing $g = g_p g_{p'}$ as a product of its $p$- and $p'$- parts gives that $g_{p'} \in C_G(x)$ and so
$x^{g_p} = x^{g} = xz$. Thus $\langle x,z, g_p \rangle$ is a $p$-group and so abelian, a contradiction.$\Box$\\

\section{Proof of Main Result}

Suppose $G$ is a finite group with either ${\rm cs}^*(G) = \{2^{\alpha}, 2^{\alpha +1}, 2^{\alpha}3\}$ or
${\rm cs}^*(G) = \{2^{\alpha}, 2^{\alpha+1}, 2^{\alpha}3, 2^{\alpha +2}\}$ and $\alpha \geq 1$. As mentioned in the previous section we
can assume $G$ is a $\{2,3\}$-group. Also consideration of  ${\rm cs}^*(G)$ makes it clear that $G$ is not nilpotent. 
We let $P$ denote a Sylow $3$-subgroup of $G$.
Recall a group $G$ is a $q$-Baer group if $q$ is a prime dividing the order of $G$ and all $q$-elements have prime power
index \cite{rachel}.\\[1ex]

{\noindent}{\bf Step 1} {\it Let $\bar{G}=G/O_2(G)Z(G)$. Then $\bar{G}$ has a normal elementary abelian Sylow 3-subgroup
$\bar{P}$ and  $|\bar{G}/\bar{P}| \leq 2$.}

Note $PZ(G)/Z(G)$ is elementary abelian by Lemma \ref{central}. Thus  
the group $\bar{G}$ has an elementary 
abelian Sylow $3$-subgroup where every $3$-element has index a power of $2$. 
Thus $\bar{G}$ is a $3$-Baer group, and as $O_2(\bar{G})=1$ so the Sylow $3$-subgroup of $\bar{G}$
is normal \cite{rachel}. Note that any 2-element of 2-power index in $G$ lies in 
$O_2(G)$, by Wielandt \cite[Lemma 6]{baer}. 
Thus $\bar{G}/\bar{P}$ is a 2-group acting on an elementary abelian $3$-group so that 
every $2$-element has index $3$.  If $t$ is any element acting on $\bar{P}$ we get 
$\bar{P}=[\bar{P},t]\oplus C_{\bar{P}}(t)$ where $|[\bar{P},t]|=3$. So $t$ inverts $[\bar{P},t]$ and has 
order $2$. Hence the determinant of the linear transformation induced by
$t$ is $-1$. This holds for every element in $\bar{G}/\bar{P}$ and thus $|\bar{G}/\bar{P}| \leq 2$.\\[1ex]

{\noindent}{\bf Step 2} {\it Suppose $x$ is a 3-element and $t$ is a 2-element with $[x,t]=1$ and 
$C_{O_2(G)}(x) \geq C_{O_2(G)}(t)$, then $x$ centralises $O_2(G)$.}

This is a direct consequence of  Thompson's $P \times Q$ lemma \cite[8.2.8]{kurzweil}.\\[1ex]

{\noindent}{\bf Step 3} {\it We can assume $|P|=3$.}

If all 2-elements have 2-power index then the Sylow 2-subgroup of $G$ is a direct factor 
\cite[Lemma 3]{rachel}, which
is clearly false. Choose $t$ a 2-element of index $2^{\alpha}3$. If $O_2(G)$ is not a Sylow 2-subgroup of $G$ choose 
$t \not\in O_2(G)$ (this is possible by \cite[Lemma 6]{baer}) otherwise choose $t \in O_2(G)$. Suppose 
$x$ is a non-central $3$-element in $C_G(t)$. Then, as $C_G(xt) = C_G(x) \cap C_G(t)$, it follows that
$xt$ has index $2^{\alpha}3$ and 
$C_G(x)\geq C_G(t)$. Applying Step 2 gives that 
 $x$  centralises $O_2(G)$. Also $x$ centralises $t$ and so, by Step 1,  
centralises a Sylow 2-subgroup and thus has index prime to 2, which is false.
So $C_G(t)$ contains no non-central $3$-elements. If $P_0$ is the
Sylow $3$-subgroup of $C_G(t)$ it is central in $G$. Let $P$ be a Sylow
$3$-subgroup of $G$ containing $P_0$ then $|P:P_0|=3$. Thus $P$ is abelian and applying Lemma \ref{centre} gives the result.
\\[1ex]

{\noindent}{\bf Step 4} {\it If $P$ is normal then ${\rm cs}^*(G) = \{2,4,6\}$.}

We show that $P$ normal implies $\alpha = 1$. The result then follows by the analysis in \cite{bianchi} and \cite{glasby}.
Note if $P$ is normal then $O_2(G)$ has index 2 in $G$ otherwise $G$ is nilpotent which is false.
As $P$ is abelian it follows that any 3-element of $G$ has index 2, i.e. $\alpha =1$ as claimed.\\[1ex]

{\noindent}{\bf Step 5} {\it Suppose
$P$ is not normal.
  Then the 2-elements
are either central or have index $2^{\alpha}$ or $2^{\alpha}3$, and these all occur. Furthermore the 3-elements have index $2^{\alpha +1}$.}

Clearly any element central in a Sylow 2-subgroup is central so we have central 2-elements. 
If all 2-elements have 2-power index then the Sylow 2-subgroup is a direct factor \cite[Lemma 3]{rachel}, 
which is clearly
false, so there exist 2-elements of index $2^{\alpha}3$. Also as any element of index $2^{\alpha}$ lies in the
Fitting subgroup  of $G$ \cite[Proposition 3.1]{jabara}, we have 2-elements of index $2^{\alpha}$. Note 3-elements do
not lie in the Fitting subgroup so do not have index $2^{\alpha}$. Furthermore, as $P$ is abelian 3-elements do not have index 
$2^{\alpha}3$. 

In the rank 4 case suppose there exists a 2-element of index $2^{\alpha+2}$. However then $G$ has a normal 2-complement \cite[Theorem 1]{alan}, contradicting our hypothesis.

 We need to consider the case that $t$ is a 2-element
of index $2^{\alpha +1}$. Let $x$ be a 3-element in $C_G(t)$. By the first paragraph the index of  $x$ is either $2^{\alpha +1}$ or
$2^{\alpha + 2}$.
 Now $C_G(xt) = C_G(x) \cap C_G(t)$ so
$xt$ has index $2^{\alpha + 1}$ or $2^{\alpha + 2}$. If $xt$ has index $2^{\alpha+1}$ then $C_G(xt) = C_G(x) = C_G(t)$
and $x$ centralises $O_2(G)$ by Step 2, giving $x$ has index 2, a contradiction. Suppose $xt$ has index $2^{\alpha+2}$.
 Then $C_G(xt)$ is normal in $C_G(t)$. 
 Choose $g \in C_G(t) \setminus C_G(xt)$ and let $P = \langle x \rangle \leq C_G(xt)$. Then $g$ acts on $P$ and
so inverts $x$. If $g \in O_2(G)$ then $[g,x] \in O_2(G)$ forces $x \in O_2(G)$, a contradiction. Thus all elements in $C_G(t) \cap O_2(G)$
centralise $x$. Applying Step 2 gives that 
$x$ centralises $O_2(G)$, again a contradiction.

As a Sylow 3-subgroup $P$ is cyclic of order $3$ it follows that all 3-elements have the same index. 
As mentioned previously 3-elements do not have index
$2^{\alpha}$ or $2^{\alpha}3$.  Thus, in the conjugate rank 3 case, the
result follows. For the conjugate rank 4 case suppose the 3-elements have
index $2^{\alpha+2}$. But then any element of mixed order also has this index and there are no elements of index
$2^{\alpha +1}$.\\[1ex]

{\noindent}{\bf Step 6} {\it Suppose ${\rm cs}^*(G) = \{ 2^{\alpha}, 2^{\alpha +1}, 2^{\alpha}3 \}$ then $\alpha = 1$.}

Let $P= \langle x \rangle$, if $P$ is normal the result holds by Step 4, so assume $P$ not normal. Thus $x$ has index
$2^{\alpha +1}$ by Step 5.
Let $t$ be a noncentral 2-element in $C_G(x)$, then
$|t^G| = 2^{\alpha}$ by Step 5. Now $C_G(xt) = C_G(x) \cap C_G(t) \geq P$ and thus $|(xt)^G| = 2^{\alpha +1}$ and, in
particular, $C_G(x) = C_G(xt)$ has index 2, and so is normal, in $C_G(t)$. Now we proceed as in Step 5 and show that all
elements in $C_G(t) \cap O_2(G)$ centralise $x$, then applying Step 2 gives that $x$
centralises $O_2(G)$, this contradicts $|x^G| = 2^{\alpha +1}$.\\[1ex]

Thus the first part of our Theorem is proved.
 From now on we suppose $G$ is a finite group with
${\rm cs}^*(G) = \{2^{\alpha}, 2^{\alpha +1}, 2^{\alpha}3, 2^{\alpha + 2} \}$. By Step 4 we know $P$ is not normal in $G$.\\[1ex]

{\noindent}{\bf Step 7} {\it The number of Sylow 3-subgroups of $G$ is $2^{\alpha +1}$ if $\alpha$ is odd and $2^{\alpha}$
if $\alpha$ is even. Furthermore, $\alpha$ is odd if and only if the Sylow 2-subgroup of $G$ is normal.}

Suppose $P = \langle x \rangle$, then $C_G(P) = C_G(x)$. Further, $|N_G(P): C_G(P)|=1$ or 2 since $|P| =3$. The number of
Sylow 3-subgroups, denoted $n_3$, satisfies $n_3 \equiv 1 \bmod 3$ and $|G:C_G(x)| = 2^{\alpha +1}$ by Step 5, hence
$n_3 = 2^{\alpha +1}$ if $\alpha$ is odd and $n_3 = 2^{\alpha}$ if $\alpha$ is even. Furthermore, if $|N_G(P):C_G(P)|=1$
then $G$ has a normal 3-complement by Burnside's $p$-complement theorem \cite[Theorem 7.2.1]{kurzweil}, 
so a normal Sylow 2-subgroup. Finally note that  if $t \in N_G(P) \setminus C_G(P)$ then $t \not\in O_2(G)$, so if $G$
has a normal Sylow 2-subgroup then $\alpha$ is odd.\\[1ex]

{\noindent}{\bf Step 8} {\it Suppose $G$ has a normal Sylow 2-subgroup,  then $\alpha =1$.}

As $G$ is not nilpotent, by Step 5 it
follows that the 2-elements are central or have index $2^{\alpha}$ or $2^{\alpha}3$. Let $T$ denote the unique Sylow
2-subgroup of $G$ and suppose $t$ is a non-central 
2-element then $|C_G(t): C_T(t)| = |TC_G(t): T|$ which is equal to 1 or 3. Further,
$$|G:C_G(t)||C_G(t):C_T(t)|= |G:T||T:C_T(t)|=3|T:C_T(t)|,$$ and so 
$|T:  C_T(t)| = 2^{\alpha}$. Thus $T/Z(T)$ has exponent 2, 
by \cite[Corollary 2.2]{ishi}, and is therefore
abelian. Also $Z(T) = Z(G)$. 
Let $P=\langle x \rangle$ be a Sylow 3-subgroup of $G$. Then $C_T(P)$ contains $Z(T)$ so is normalised by $T$ and $P$
so $C_T(P)$ is normal. But then $C_T(P)$ centralises all Sylow 3-subgroups.
Now choose $y \in C_T(P)$  such that $xy$ has index $2^{\alpha+2}$, such a $y$ exists by Step 5.
So, $|C_G(y): C_G(x) \cap C_G(y)| = 4$ and thus $y$ commutes with exactly 4 Sylow 3-subgroups. However $y \in C_T(P)$
so $y$ centralises all Sylow 3-subgroups. By the previous step, $G$ has $2^{\alpha+1}$ Sylow 3-subgroups. So
 $2^{\alpha +1} = 4$
and $\alpha =1$.\\[1ex]

This completes the proof of our Theorem.

\noindent Rachel D. Camina: Fitzwilliam College, Cambridge, CB3 0DG, UK.\\
rdc26@cam.ac.uk\\[1ex]
Alan R. Camina: School of Mathematics, UEA, Norwich, NR4 7TJ, UK.\\
A.Camina@uea.ac.uk

\end{document}